\documentclass[a4paper, 10pt, conference]{ieeeconf}      

\IEEEoverridecommandlockouts                              

\usepackage{amsmath,amsfonts,amssymb}
\usepackage{graphicx}
\usepackage{cases}
\usepackage{algorithm}
\usepackage{algorithmicx}
\usepackage{algpseudocode}
\algrenewcommand\algorithmicrequire{\textbf{Input:}}
\algrenewcommand\algorithmicensure{\textbf{Output:}}
\usepackage{siunitx}
\usepackage{wrapfig}
\usepackage{subcaption}
\usepackage{fancyhdr}
\fancypagestyle{fancyheading}{
  \fancyfoot{}
}

\newtheorem{thm}{Theorem}[section]

\newtheorem{prop}[thm]{Proposition}
\newtheorem{cor}[thm]{Corollary}

\newtheorem{defn}{Definition}[section]

\newcommand{\norm}[1]{\left\lVert#1\right\rVert}

\title{\LARGE \bf
Stabilizing Traffic via Autonomous Vehicles: A Continuum Mean Field Game Approach
}
\author{Kuang Huang$^{1}$, Xuan Di$^{2,3}$, Qiang Du$^{1,3}$ and Xi Chen$^{4}$
\thanks{$^{1}$Department of Applied Physics and Applied Mathematics, Columbia University}
\thanks{$^{2}$Department of Civil Engineering and Engineering Mechanics, Columbia University}
\thanks{$^{3}$Data Science Institute, Columbia University}
\thanks{$^{4}$Department of Computer Science, Columbia University}
}

\begin{document}

\maketitle
\setlength{\topmargin}{-48pt}
\setlength{\headsep}{40pt}
\thispagestyle{fancyheading}
\pagestyle{fancyheading}

\begin{abstract}
	This paper presents scalable traffic stability analysis for both pure autonomous vehicle (AV) traffic and mixed traffic based on continuum traffic flow models. Human vehicles are modeled by a non-equilibrium traffic flow model, i.e., Aw-Rascle-Zhang (ARZ), which is unstable. AVs are modeled by the mean field game which assumes AVs are rational agents with anticipation capacities. It is shown from linear stability analysis and numerical experiments that AVs help stabilize the traffic. Further, we quantify the impact of AV's penetration rate and controller design on the traffic stability. The results may provide insights for AV manufacturers and city planners.

\end{abstract}

\section{Introduction}
\label{sec:intro}
	
Autonomous vehicles (AVs) are believed to be the foundation of the next-decade transportation system 
and are expected to improve the traffic flow that is presently dominated by human vehicles (HVs). 
Modeling AV's driving behavior and quantifying different penetration rates of AVs' impact on the traffic is of great significance.

This paper focuses on traffic stability, which is one of the most substantial traffic features. Traffic stability refers to a traffic system's asymptotic stability around uniform flows. HV traffic is observed to be an unstable system in which a small perturbation (caused by driving errors or delays) to the uniform flow will grow up with time and develop traffic congestion. By removing human errors, it is expected that AVs
will help stabilize the traffic system. A field experiment \cite{stern2018dissipation} showed that one AV is able to stabilize the traffic with approximately twenty vehicles on a ring road.

AV's capability of stabilizing traffic is validated using microscopic models for mixed AV-HV traffic. In microscopic models, the traffic system is described by ordinary differential equations.
One can carry out standard linear stability analysis to characterize such a traffic system's stability, built upon 
connected cruise controllers \cite{jin2018connected} or generic car-following models \cite{cui2017stabilizing,wu2018stabilizing}. \cite{jin2018connected,cui2017stabilizing} considered only one AV with multiple HVs; \cite{wu2018stabilizing} studied multiple AVs and multiple HVs but only focused on the head-to-tail stability. 
However, in the general case, the mixed traffic stability analysis relies on the topology of mixed vehicles and the vehicle-to-vehicle communication network,
which suffers from scalability issues. 

One alternative approach to address the scalability issues is the PDE approximation \cite{zheng2016stability}. 
This approach suggests to study the stability of continuum traffic flow models which are the limits of microscopic models. The approach is well suited for the mixed traffic since one needs only to concern about the density distributions of different classes. 

In continuum traffic flow models, 
the traffic system is described by partial differential equations (PDEs) of traffic density and velocity. For single class traffic, the Lighthill-Whitham-Richards (LWR) model \cite{lighthill1955kinematic} is the most extensively used continuum model.
As a generalization, the multiclass LWR is widely used to model the interaction between two types of vehicles. 
\cite{levin2016multiclass} is among a few studies that applied multi-class LWR models to AV-HV mixed traffic and proposed networked traffic controls in the presence of AVs. 
Based on gas-kinetic theory, \cite{ngoduy2009continuum} proposed a multiclass macroscopic model to capture the effect of communication and information sharing on traffic flow
and analyzed the model's stability with respect to the connected vehicle's penetration rate; 
\cite{porfyri2015stability,delis2016simulation,delis2018macroscopic} modeled the macroscopic traffic flow of mixed Adaptive Cruise Control (ACC) and Cooperative Adaptive Cruise Control (CACC) vehicles and analyzed how the ACC vehicle's penetration rate influences the traffic stability.

This paper models AVs using the \emph{mean field game} following the authors' work \cite{huang2019game}. In this framework, AVs are assumed to be rational, utility-optimizing agents with anticipation capabilities and play a non-cooperative game by selecting their driving speeds. 
AVs' utility-optimizing and anticipation behaviors are distinctive characteristics from the aforementioned continuum models. By extending \cite{huang2019game}, this paper aims to build continuum traffic flow models for both pure AV traffic and mixed AV-HV traffic based on mean field games and analyze the models' traffic stability.

The remainder of the paper is organized as follows. Section \ref{sec:pre} provides an overview of the mean field game and the Aw-Rascle-Zhang model, used for modeling AVs and HVs, respectively. Section \ref{sec:form} formulates models for both pure AV traffic and mixed AV-HV traffic. Based on the proposed models, Section \ref{sec:av_analysis} shows the linear stability analysis for the pure AV traffic and Section \ref{sec:mixed_analysis} demonstrates the mixed traffic's stability through numerical experiments.
\section{Preliminaries}
\label{sec:pre}

\subsection{Mean Field Game}
\label{sec:pre_mfg}

Mean field game (MFG) is a game-theoretic framework to model complex multi-agent dynamic systems \cite{lasry2007mean}. 
In the MFG framework, a population of $N$ rational utility-optimizing agents are modeled by a dynamic system. The agents interact with each other through their utilities. Assuming those agents optimize their utilities in a non-cooperative way, they form a \emph{differential game}.

Exact Nash equilibria to the differential game are generally hard to solve when $N$ is large. Alternatively, MFG considers the continuum problem as $N\to\infty$. By exploiting the ``smoothing" effect of a large number of interacting individuals, MFG assumes that each agent only responds to and contributes to the density distribution of the whole population. Then the equilibria are characterized by a set of two PDEs: a backward Hamilton-Jacobi-Bellman (HJB) equation describing a generic agent's optimal control provided the density distribution and a forward Fokker-Planck equation describing the population's density evolution provided individual controls.


In this paper we shall formulate AV traffic as a mean field game. AVs are modeled as rational agents with predefined driving costs. 
Their density distribution is exactly the traffic density and the Fokker-Planck equation is the same as the continuity equation (CE) that is widely used in continuum traffic flow models.  

\subsection{Aw-Rascle-Zhang Model}
\label{sec:pre_arz}
	
The following Aw-Rascle-Zhang (ARZ) model:
\begin{align}
	\mbox{(CE)}\quad \rho_t+(\rho u)_x&=0,\label{eq:arz1}\\
	\mbox{(ME)}\quad [u+h(\rho)]_t+u[u+h(\rho)]_x&=\frac{1}{\tau}[U(\rho)-u],\label{eq:arz2}
\end{align}
is a non-equilibrium continuum traffic flow model describing human driving behaviors \cite{aw2000resurrection,zhang2002non}, where,\\
$\rho(x,t)$, $u(x,t)$: the traffic density and speed;\\
$U(\cdot)$: the desired speed function;\\
$h(\cdot)$: the hesitation function that is an increasing function of the density;\\
$\tau$: the relaxation time quantifying how fast drivers adapt their current speeds to desired speeds.

Equation \eqref{eq:arz1} is the continuity equation describing the flow conservation and \eqref{eq:arz2} is a momentum equation (ME) prescribing human driver's dynamic behavior.

The ARZ model is able to predict important human driving features such as stop-and-go waves and traffic instability \cite{Seibold2012}. Traffic stability is defined around uniform flows. In continuum models, uniform flows are described by constant solutions $\rho(x,t)\equiv\bar{\rho}$, $u(x,t)\equiv\bar{u}$. The constant solutions of the ARZ model is given by $\bar{u}=U(\bar{\rho})$. Then the traffic stability for the ARZ model is defined as follows:

\begin{defn}\label{def:stab1}
	The ARZ model \eqref{eq:arz1}\eqref{eq:arz2} is stable around the uniform flow $(\bar{\rho},\bar{u})$ where $\bar{u}=U(\bar{\rho})$ if for any $\varepsilon>0$, there exists $\delta>0$ such that for any solution $\rho(x,t),u(x,t)$ to the system:
	\begin{align}
		\sup_{0\leq t<\infty}\left\{\norm{\rho(\cdot,t)-\bar{\rho}}+\norm{u(\cdot,t)-\bar{u}}\right\}\leq \varepsilon,
	\end{align}
	whenever $\norm{\rho(\cdot,0)-\bar{\rho}}+\norm{u(\cdot,0)-\bar{u}}\leq \delta$.
	Here $\norm{\cdot}$ is a given norm. The system is linearly stable if its linearized system at $(\bar{\rho},\bar{u})$ is stable around the zero solution.
\end{defn}

 The ARZ model has a simple linear stability criterion \cite{Seibold2012}:
\begin{thm}
	The ARZ model \eqref{eq:arz1}\eqref{eq:arz2} is linearly stable around the uniform flow $(\bar{\rho},\bar{u})$ where $\bar{u}=U(\bar{\rho})$ if and only if $h'(\bar{\rho})>-U'(\bar{\rho})$.
\end{thm}

Because of its capability of producing traffic instability, we shall use the ARZ model \eqref{eq:arz1}\eqref{eq:arz2} to characterize HV's driving behavior.
\section{Model Formulation}
\label{sec:form}

\subsection{Pure AV Traffic: Mean Field Game}
In this section, we will build a pure AV continuum traffic flow model based on a mean field game following \cite{huang2019game}. 

Assume that a large population of homogeneous AVs are driving on a closed highway without any entrance nor exit. Those AVs anticipate others' behaviors and the evolution of the traffic density $\rho(x,t)$ on a predefined time horizon $[0,T]$. AVs control their speeds and aim to minimize their driving costs on the horizon $[0,T]$. 
Then the AVs' optimal cost $V(x,t)$ and optimal velocity field $u(x,t)$ can be described by a set of HJB equations \cite{huang2019game}:
\begin{align}
	\mbox{(HJB)}\quad &V_t+uV_x+f(u,\rho)=0,\label{eq:hjb1}\\
	&u=\text{argmin}_{\alpha} \{\alpha V_x+f(\alpha,\rho)\},\label{eq:hjb2}	
\end{align}
where $f(\cdot,\cdot)$ is the \emph{cost function} \cite{huang2019game}. 

When all AVs follow their optimal velocity controls, the system's density evolution is described by the continuity equation:
\begin{align}
	\mbox{(CE)}\quad \rho_t+(\rho u)_x=0.\label{eq:cemfg}
\end{align}

The mean field game is described by the coupled system \eqref{eq:hjb1}\eqref{eq:hjb2}\eqref{eq:cemfg}. 

\begin{itemize}
	\item The initial condition for the forward continuity equation \eqref{eq:cemfg} is given by the initial density $\rho(x,0)=\rho_0(x)$. 
	\item The terminal condition for the backward HJB equations \eqref{eq:hjb1}\eqref{eq:hjb2} is given by the terminal cost $V(x,T)=V_T(x)$. We will always set $V_T(x)=0$ meaning that the cars have no preference on their destinations.
	\item The choice of the spatial boundary condition depends on the traffic scenario. In this paper we assume that the highway is a ring road of fixed length $L$ and specify the periodic boundary condition $\rho(0,t)=\rho(L,t)$, $V(0,t)=V(L,t)$.
\end{itemize}

The cost function represents certain driving objectives. The choice of the cost function determines AV's driving behavior. In this paper we shall follow \cite{huang2019game} and take the following cost function:
\begin{align}
	f(u,\rho)=\underbrace{\frac12\left(\frac{u}{u_{\text{max}}}\right)^2}_{\text{kinetic energy}}-\underbrace{\frac{u}{u_{\text{max}}}}_{\text{efficiency}}+\underbrace{\frac{u\rho}{u_{\text{max}}\rho_{\text{jam}}}}_{\text{safety}},\label{eq:costfct1}
\end{align}
where,\\
$u_\text{max}$ and $\rho_\text{jam}$ are the free flow speed and the jam density;\\
$\frac12\left(u/u_{\text{max}}\right)^2$ models the car's kinetic energy;\\
$-u/u_{\text{max}}$ models the car's efficiency, minimizing this term means that the car should drive as fast as possible;\\
$u\rho/u_{\text{max}}\rho_{\text{jam}}$ models the safety, it is a penalty term that restricts the car's speed in traffic congestion;

The MFG system corresponding to the cost function \eqref{eq:costfct1} is \cite{huang2019game}:
\begin{subequations}
\begin{numcases}{}
	\rho_t+(\rho u)_x=0,\label{eq:mfg_cont}\\
	V_t+uV_x+\frac12\left(\frac{u}{u_{\text{max}}}\right)^2-\frac{u}{u_{\text{max}}}+\frac{u\rho}{u_{\text{max}}\rho_{\text{jam}}}=0,\quad\quad\label{eq:mfg_hjb1}\\
	u=g_{[0,u_\text{max}]}\left(u_{\text{max}}\left(1-\frac{\rho}{\rho_{\text{jam}}}-u_{\text{max}}V_x\right)\right),\label{eq:mfg_hjb2}
\end{numcases}
\end{subequations}
where $g_{[0,u_\text{max}]}(u)=\max\{\min\{u,u_\text{max}\},0\}$ is a cut-off function which ensures the cars' speeds satisfy the constraint $0\leq u\leq u_\text{max}$.

\cite{huang2019game} provides theoretical and numerical analysis on the MFG system \eqref{eq:mfg_cont}\eqref{eq:mfg_hjb1}\eqref{eq:mfg_hjb2}. 

The uniform flows of the MFG system \eqref{eq:mfg_cont}\eqref{eq:mfg_hjb1}\eqref{eq:mfg_hjb2} are given by $\bar{u}=u_\text{max}\left(1-\frac{\bar{\rho}}{\rho_\text{jam}}\right)$.
Note that Definition \ref{def:stab1} does not apply to the MFG system since the system is defined and solved on a fixed time horizon $[0,T]$. In this case, we define traffic stability as follows:

\begin{defn}\label{def:stab2}
	The MFG system \eqref{eq:mfg_cont}\eqref{eq:mfg_hjb1}\eqref{eq:mfg_hjb2} is stable around the uniform flow $(\bar{\rho},\bar{u})$ where $\bar{u}=u_\text{max}(1-\bar{\rho}/\rho_\text{jam})$ if for any $\varepsilon>0$, there exists $\delta>0$ such that for any $T>0$ and for any solution $\rho^{(T)}(x,t),u^{(T)}(x,t)$ to the system with $V_T(x)=0$ on the time horizon $[0,T]$:
	\begin{align}
		\sup_{0\leq t\leq T}\left\{\norm{\rho^{(T)}(\cdot,t)-\bar{\rho}}+\norm{u^{(T)}(\cdot,t)-\bar{u}}\right\}\leq \varepsilon,
	\end{align}
	whenever $\norm{\rho(\cdot,0)-\bar{\rho}}\leq \delta$.
	The system is linearly stable if its linearized system at $(\bar{\rho},\bar{u})$ is stable around the zero solution.
\end{defn}


\subsection{Mixed Traffic: Coupled MFG-ARZ System}

This section aims to develop a continuum mixed AV-HV traffic flow model. We denote $\rho^{\text{AV}}(x,t)$ the AV density, $\rho^{\text{HV}}(x,t)$ the HV density and
\begin{align}
	\rho^{\text{TOT}}(x,t)=\rho^{\text{AV}}(x,t)+\rho^{\text{HV}}(x,t),
\end{align}
the total density. Denote $u^{\text{AV}}(x,t)$ and $u^{\text{HV}}(x,t)$ the velocities of AVs and HVs, respectively.

We model HVs by the ARZ model and AVs by the MFG, respectively. The next step is to model the interactions between AVs and HVs. The interactions include the flow interaction and the dynamic interaction. 

\emph{Flow interaction}. The flow interaction relates to how the multiclass flows are computed and assigned. We follow the framework from \cite{fan2015heterogeneous} and suppose that the multiclass flows are described by the following continuity equations for both AVs and HVs: 
\begin{align}
	\mbox{(CE-AV)}\quad &\rho^{\text{AV}}_t+(\rho^{\text{AV}} u^{\text{AV}})_x=0,\\
	\mbox{(CE-HV)}\quad &\rho^{\text{HV}}_t+(\rho^{\text{HV}} u^{\text{HV}})_x=0.
\end{align}

\emph{Dynamic interaction}. Each of the velocities $u^{\text{AV}}$ and $u^{\text{HV}}$ should depend on both AV density $\rho^{\text{AV}}$ and HV density $\rho^{\text{HV}}$. The way of defining the velocities over multiclass densities characterizes the dynamic interaction. \cite{fan2015heterogeneous} summarized some possible formulations of the dynamic interaction.

In this paper we model an asymmetric dynamic interaction between AVs and HVs by introducing multiclass densities into the HJB equations and the momentum equation of the system. For HVs, we assume that HVs only observe the total density $\rho^\text{TOT}$ to adapt their speeds. The momentum equation \eqref{eq:arz2} in the ARZ model then becomes:
\begin{align}
	\left[u^{\text{HV}}+h(\rho^{\text{TOT}})]_t+u^{\text{HV}}[u^{\text{HV}}+h(\rho^{\text{TOT}})\right]_x=\notag\\
	\frac{1}{\tau}\left[U(\rho^{\text{TOT}})-u^{\text{HV}}\right].
\end{align}
We take the Greenshields desired speed function $U(\rho)=u_{\text{max}}\left(1-\rho/\rho_{\text{jam}}\right)$.

For AVs, we assume that AVs observe both AV and HV densities. We model AVs' reaction to multiclass densities by introducing an extra term into the AV's cost function. The AV's modified cost function for mixed traffic is:
\begin{align}
	f(u^{\text{AV}},\rho^{\text{AV}},\rho^{\text{HV}})=&\underbrace{\frac12\left(\frac{u^{\text{AV}}}{u_{\text{max}}}\right)^2}_{\text{kinetic energy}}-\underbrace{\frac{u^{\text{AV}}}{u_{\text{max}}}}_{\text{efficiency}}\notag\\
	+&\underbrace{\frac{u^{\text{AV}}\rho^{\text{TOT}}}{u_{\text{max}}\rho_{\text{jam}}}+\beta\frac{\rho^{\text{HV}}}{\rho_\text{jam}}}_{\text{safety}}\label{eq:costcc},
\end{align}
where the safety is modeled by two penalty terms: one is similar to the penalty term in \eqref{eq:costfct1} but the congestion is modeled by the total density $\rho^{\text{TOT}}$, the other quantifies HV's impact on AV's speed selection and the parameter $\beta$ represents AV's sensitivity to HV's density. From \eqref{eq:costcc} we can derive the corresponding HJB equations.


Summarizing all above, we obtain the following coupled MFG-ARZ system:

\begin{subequations}
\begin{numcases}{}
	\rho^{\text{AV}}_t+(\rho^{\text{AV}} u^{\text{AV}})_x=0,\label{eq:mixed_1}\\
	V_t+u^{\text{AV}}V_x+\frac12\left(\frac{u^{\text{AV}}}{u_\text{max}}\right)^2-\frac{u^{\text{AV}}}{u_\text{max}}+\frac{u^{\text{AV}}\rho^{\text{TOT}}}{u_\text{max}\rho_\text{jam}}\notag\\
	+\beta\frac{\rho^{\text{HV}}}{\rho_\text{jam}}=0,\quad\quad\ \label{eq:mixed_2}\\
	u^{\text{AV}}=g_{[0,u_\text{max}]}\left(u_\text{max}\left(1-\frac{\rho^{\text{TOT}}}{\rho_\text{jam}}-u_\text{max}V_x\right)\right),\quad\quad \label{eq:mixed_3}\\
	\rho^{\text{HV}}_t+(\rho^{\text{HV}} u^{\text{HV}})_x=0,\label{eq:mixed_4}\\
	\left[u^{\text{HV}}+h(\rho^{\text{TOT}})]_t+u^{\text{HV}}[u^{\text{HV}}+h(\rho^{\text{TOT}})\right]_x=\notag\\
	\frac{1}{\tau}\left[u_\text{max}\left(1-\frac{\rho^{\text{TOT}}}{\rho_\text{jam}}\right)-u^{\text{HV}}\right],\label{eq:mixed_5}\\
	\rho^{\text{TOT}}=\rho^{\text{AV}}+\rho^{\text{HV}}.\label{eq:mixed_6}
\end{numcases}
\end{subequations}

\begin{itemize}
	\item The initial conditions are given by the initial densities $\rho^\text{AV}(x,0)=\rho^\text{AV}_0(x)$, $\rho^\text{HV}(x,0)=\rho^\text{HV}_0(x)$ and the initial velocity $u^\text{HV}(x,0)=u^\text{HV}_0(x)$.
	\item The terminal condition is given by the terminal cost $V(x,T)=V_T(x)$. We will always set $V_T(x)=0$.
	\item We specify the periodic boundary conditions for all of $\rho^{\text{AV}}$, $\rho^{\text{HV}}$, $u^{\text{HV}}$ and $V$.
\end{itemize}

	The mixed traffic's uniform flows are defined as the system's constant solutions $\rho^{\text{AV}}(x,t)\equiv\bar{\rho}^{\text{AV}}$, $\rho^{\text{HV}}(x,t)\equiv\bar{\rho}^{\text{HV}}$,
	\begin{align}
	 	\rho^{\text{TOT}}(x,t)\equiv\bar{\rho}^{\text{TOT}}=\bar{\rho}^{\text{AV}}+\bar{\rho}^{\text{HV}},
	\end{align}
	and 
	\begin{align}
		u^{\text{AV}}(x,t)\equiv u^{\text{HV}}(x,t)\equiv\bar{u}=u_\text{max}\left(1-\frac{\bar{\rho}^\text{TOT}}{\rho_\text{jam}}\right).\label{eq:mix_con_sol}
	\end{align}
	Since AVs are modeled by a mean field game, the mixed traffic system (\ref{eq:mixed_1}-\ref{eq:mixed_6}) is defined and solved on a predefined time horizon $[0,T]$. Similar to Definition \ref{def:stab2}, the mixed traffic system's stability is defined as:
\begin{defn}\label{def:stab3}
	The system (\ref{eq:mixed_1}-\ref{eq:mixed_6}) is stable around the uniform flow $(\bar{\rho}^{\text{AV}},\bar{\rho}^{\text{HV}},\bar{u})$ which satisfies \eqref{eq:mix_con_sol} if for any $\varepsilon>0$, there exists $\delta>0$ such that for any $T>0$ and for any solution $\rho^{\text{AV},(T)}(x,t)$, $u^{\text{AV},(T)}(x,t)$, $\rho^{\text{HV},(T)}(x,t)$, $u^{\text{HV},(T)}(x,t)$ to the system with $V_T(x)=0$ on the time horizon $[0,T]$:
	\begin{align}
		\sup_{0\leq t\leq T}\sum_{i=\text{AV},\text{HV}}\norm{\rho^{i,(T)}(\cdot,t)-\bar{\rho}^i}+\norm{u^{i,(T)}(\cdot,t)-\bar{u}}\leq \varepsilon,\label{eq:stab1}
	\end{align}
	whenever
	\begin{align}
		\sum_{i=\text{AV},\text{HV}}\norm{\rho^{i,(T)}(\cdot,0)-\bar{\rho}^i}+\norm{u^{\text{HV},(T)}(\cdot,0)-\bar{u}}\leq \delta.\label{eq:stab2}
	\end{align}
	The system is linearly stable if its linearized system at $(\bar{\rho}^{\text{AV}},\bar{\rho}^{\text{HV}},\bar{u})$ is stable around the zero solution.
\end{defn}

\section{Pure AV Traffic: Linear Stability Analysis}
\label{sec:av_analysis}
In this section we will carry out the standard linear stability analysis for the MFG system \eqref{eq:mfg_cont}\eqref{eq:mfg_hjb1}\eqref{eq:mfg_hjb2}. 

By scaling to dimensionless quantities we assume $u_{\text{max}}=1$ and $\rho_{\text{jam}}=1$. In addition we remove the speed constraint $0\leq u\leq u_\text{max}$ since the existence of the constraint does not change the system's stability when $0<\bar{u}<u_\text{max}$. Then we eliminate $V$ from the system \eqref{eq:mfg_cont}\eqref{eq:mfg_hjb1}\eqref{eq:mfg_hjb2} and obtain a simpler system of $\rho$ and $u$:
\begin{align}
\begin{cases}
	\rho_t+(\rho u)_x=0,\\
	u_t+uu_x-(\rho u)_x=0.
\end{cases}\label{mfg_clean}
\end{align}

Fix a uniform flow $(\bar{\rho},\bar{u})$ where $\bar{u}=1-\bar{\rho}$. Suppose that the system \eqref{mfg_clean} has the initial condition $\rho(x,0)=\bar{\rho}+\tilde{\rho}_0(x)$ and the terminal condition $V_T(x)=0$. Here $\tilde{\rho}_0(x)$ is any small perturbation.


Then we linearize the system \eqref{mfg_clean} near the uniform flow $(\bar{\rho},\bar{u})$. Suppose $\rho(x,t)=\bar{\rho}+\tilde{\rho}(x,t)$, $u(x,t)=\bar{u}+\tilde{u}(x,t)$. Note that $\bar{u}=1-\bar{\rho}$, we get the following linearized system:
\begin{align}
	\begin{cases}
	\tilde{\rho}_t+(1-\bar{\rho})\tilde{\rho}_x+\bar{\rho}\tilde{u}_x=0,\\
	\tilde{u}_t+(\bar{\rho}-1)\tilde{\rho}_x+(1-2\bar{\rho})\tilde{u}_x=0.
\end{cases}\label{mfg_lin}
\end{align}
\eqref{mfg_lin} is also a forward-backward system with the initial condition $\tilde{\rho}(x,0)=\tilde{\rho}_0(x)$ and the terminal condition $\tilde{\rho}(x,T)+\tilde{u}(x,T)=0$.

\begin{prop}\label{prop:lin}
	The linearized system \eqref{mfg_lin} is stable near the zero solution for all $0<\bar{\rho}<1$.
\end{prop}

We provide a computer-assisted proof for Proposition \ref{prop:lin} in the Appendix. The analytical proof is left for future research. As a corollary of Proposition \ref{prop:lin} we have the following results on the MFG system's stability:
\begin{cor}
	The MFG system \eqref{mfg_clean} is linearly stable around the uniform flow $(\bar{\rho},\bar{u})$ where $\bar{u}=1-\bar{\rho}$ for all $0<\bar{\rho}<1$.
\end{cor}

Our analysis shows that the proposed MFG system for AVs is always stable even if each AV only aims to optimize his own utility. Then we turn our attention to the mixed traffic and study whether the existence of AVs can stabilize the unstable HV traffic.
\section{Mixed Traffic: Numerical Experiments}
\label{sec:mixed_analysis}
In this section, we will demonstrate the stability of the mixed traffic system (\ref{eq:mixed_1}-\ref{eq:mixed_6}) by numerical experiments. We will run numerical simulations in different scenarios and check the stability in those simulations automatically with a stability criterion. Then we discuss how AVs' different penetration rates and different controller designs influence the stabilizing effect. 

\subsection{Experimental Settings} 

Take vehicles' free flow speed $u_{\text{max}}=\SI{30}{m/s}$ and the jam density $\rho_{\text{jam}}=\SI{1/7.5}{m}$. Choose the hesitation function $h(\rho)$ in the ARZ model to be:
\begin{align}
	h(\rho)=\SI{9}{m/s}\cdot\left(\frac{\rho/\rho_\text{jam}}{1-\rho/\rho_\text{jam}}\right)^{1/2},
\end{align}
which has the same form as the one used in \cite{Seibold2012}. For all of the numerical experiments, the length of the ring road $L=\SI{1}{km}$ and the length of the time horizon $T=2L/u_{\text{max}}$. 

For the system (\ref{eq:mixed_1}-\ref{eq:mixed_6}) and its arbitrary uniform flow solution $(\bar{\rho}^{\text{AV}},\bar{\rho}^{\text{HV}},\bar{u})$, the initial densities are set to be:
\begin{align}
	\rho^i_0(x)=\bar{\rho}^i+0.1\times\bar{\rho}^i\sin(2\pi x/L),
\end{align}
for $i=\text{AV},\text{HV}$ so that the initial perturbations on both AV and HV densities are sine waves whose magnitudes are 10\% of the respective uniform states. The HV's initial velocity is set to be:
\begin{align}
	u^{\text{HV}}_0(x)\equiv\bar{u}=u_\text{max}\left(1-\frac{\bar{\rho}^\text{TOT}}{\rho_\text{jam}}\right),
\end{align}
where $\bar{\rho}^\text{TOT}=\bar{\rho}^\text{AV}+\bar{\rho}^\text{HV}$ so that there is no initial perturbation on HV's velocity.
The AV's terminal cost is always set to be $V_T(x)=0$.

It is not easy to check the conditions \eqref{eq:stab1}\eqref{eq:stab2} directly. Alternatively we shall use a simplified stability criterion. Suppose $\rho^{\text{AV},(T)}(x,t)$, $u^{\text{AV},(T)}(x,t)$, $\rho^{\text{HV},(T)}(x,t)$, $u^{\text{HV},(T)}(x,t)$ is any solution to the system, we define an \emph{error function}:
\begin{align}
	E(t)=\sum_{i=\text{AV},\text{HV}}\norm{\rho^{i,(T)}(\cdot,t)-\bar{\rho}^i}+\norm{u^{i,(T)}(\cdot,t)-\bar{u}},
\end{align}
for $0\leq t\leq T$ and the system is said to be unstable if:
\begin{align}
	\max_{0\leq t\leq T}E(t)\geq 2E(0),\label{eq:crit}
\end{align}
otherwise it is said to be stable. The stability criterion \eqref{eq:crit} is checked automatically in the numerical experiments. It is validated in the experiments with no presence of AVs that the criterion \eqref{eq:crit} predicts the same stability as the ARZ model's analytical stability criterion.

\subsection{Numerical Method}
	To solve the coupled MFG-ARZ system (\ref{eq:mixed_1}-\ref{eq:mixed_6}) numerically, we apply a finite difference method (FDM) on the spatial-temporal grids. 
	We discretize the continuity equations \eqref{eq:mixed_1}\eqref{eq:mixed_4} by the Lax-Friedrichs scheme. We discretize the HJB equations \eqref{eq:mixed_2}\eqref{eq:mixed_3} of the MFG by an upwind scheme \cite{huang2019game}.  The momentum equation \eqref{eq:mixed_5} of the ARZ model is transformed into its conservative form with a relaxation term. Then we apply a hybrid scheme with an explicit Lax-Friedrichs scheme for the conservation part and an implicit Euler scheme for the relaxation part. Finally we compress all equations into a large nonlinear system and solve the system by Newton's method \cite{huang2019game}. 

\subsection{Numerical Results}

\begin{figure}[htbp]
\centering
\begin{subfigure}{0.2\textwidth}
\centering
    \includegraphics[width=\textwidth]{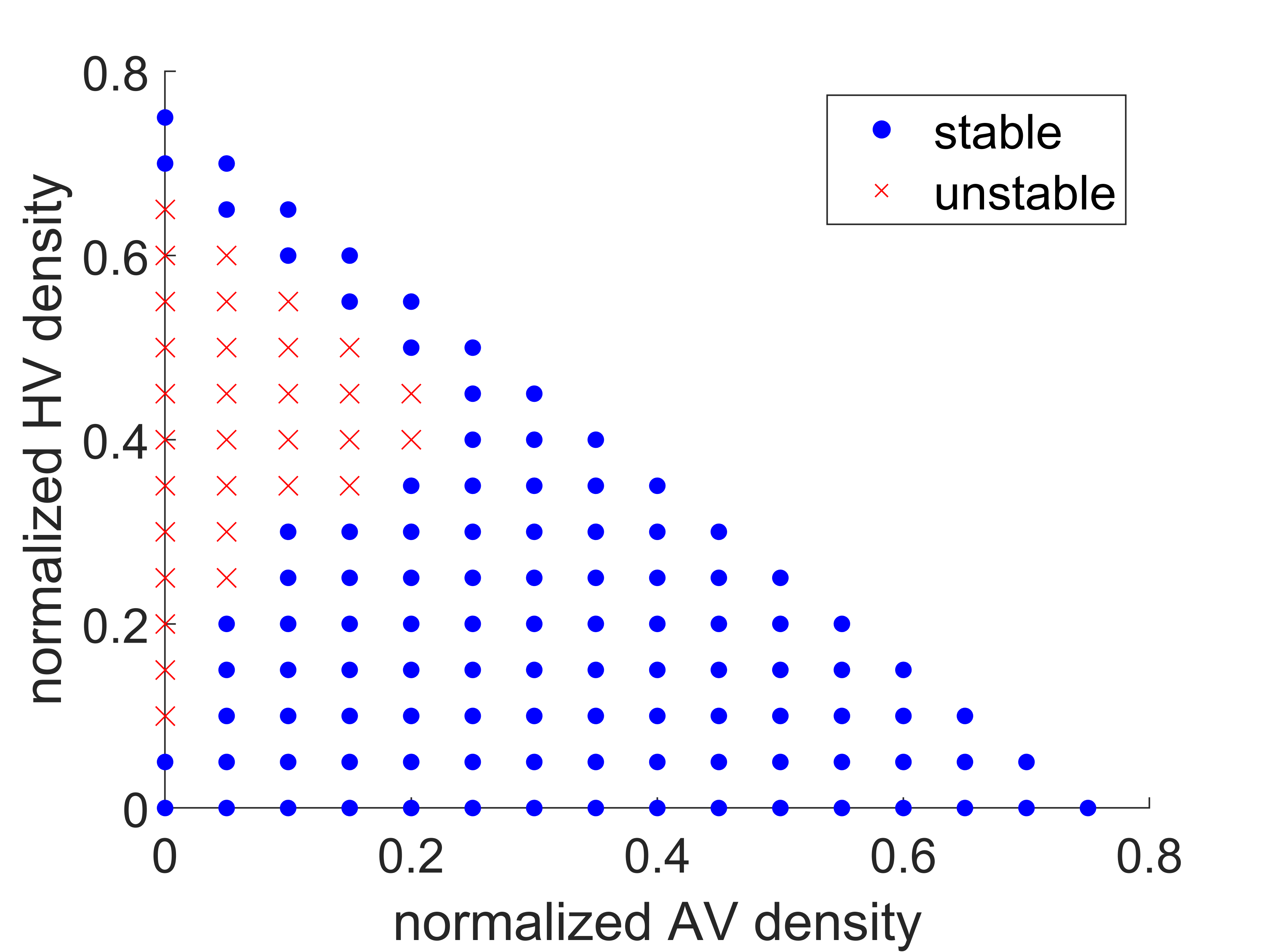}
    \caption{First group}
    \label{fig:avhv}
\end{subfigure}%
\begin{subfigure}{0.2\textwidth}
\centering
    \includegraphics[width=\textwidth]{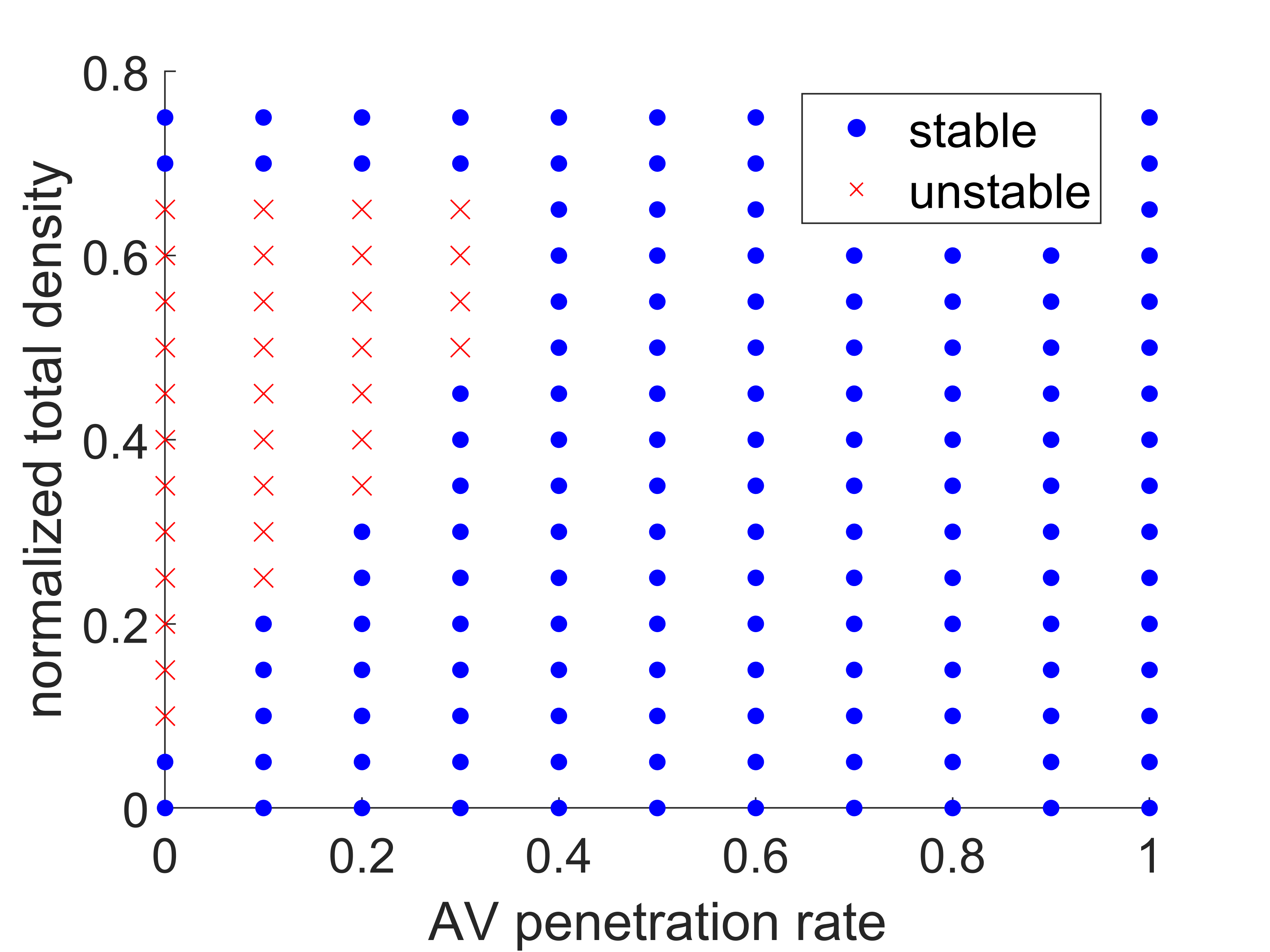}
    \caption{Second group}
    \label{fig:pc}
\end{subfigure}
\caption{Stability regions for the first and second groups of experiments}
\end{figure}

In the first group of experiments we fix $\beta=0$ and try different pairs of $\bar{\rho}^{\text{AV}}$ and $\bar{\rho}^{\text{HV}}$. We restrict the values to be under $\bar{\rho}^{\text{AV}}+\bar{\rho}^{\text{HV}}\leq0.75\rho_{\text{jam}}$ to avoid the total density exceeding the jam density. We check the system's stability from each numerical experiment and plot the results in the phase diagram between the normalized AV and HV density, see Figure \ref{fig:avhv}. We observe that when the HV density is fixed, adding AVs can stabilize the traffic. 
when the AV density is large enough, the mixed traffic is always stable.

In the second group of experiments we still keep $\beta=0$ but try different total densities $\bar{\rho}^{\text{TOT}}$ and different AV's penetration rates. Then we plot the results in the phase diagram between the AV's penetration rate and the normalized total density, see Figure \ref{fig:pc}. We observe that when the total density is fixed, traffic becomes more stable with a higher portion of AVs. In addition, the minimal AV's penetration rate to make the traffic stable increases as the total density increases. We also observe that when the AV's penetration rate is large enough, the mixed traffic is always stable.

\begin{figure}[htbp]
	\centering
	\begin{minipage}[t]{.2\textwidth}
	\centering
	\includegraphics[width=.8\textwidth]{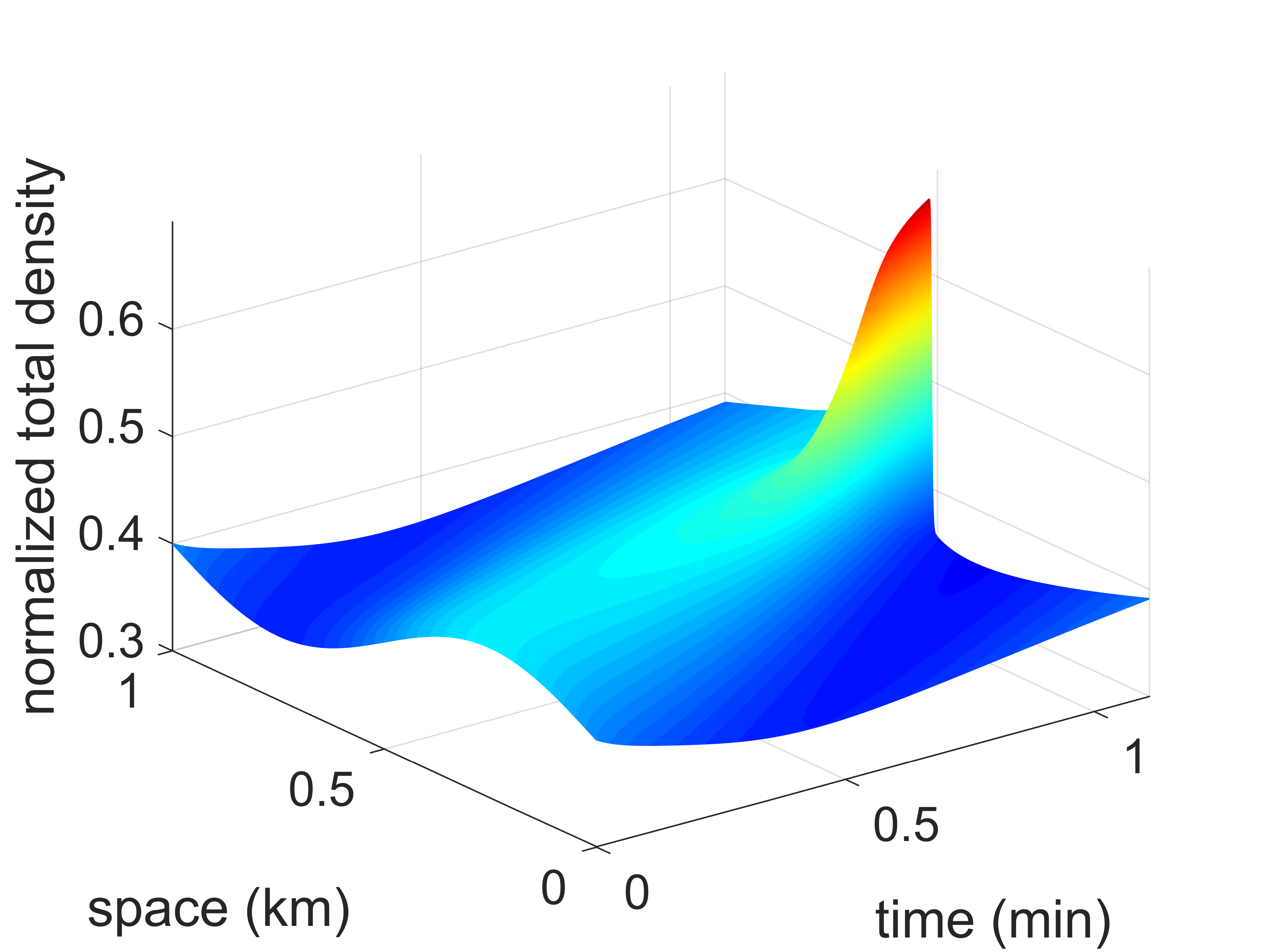}
	\end{minipage}
	\begin{minipage}[t]{.2\textwidth}
	\centering
	\includegraphics[width=.8\textwidth]{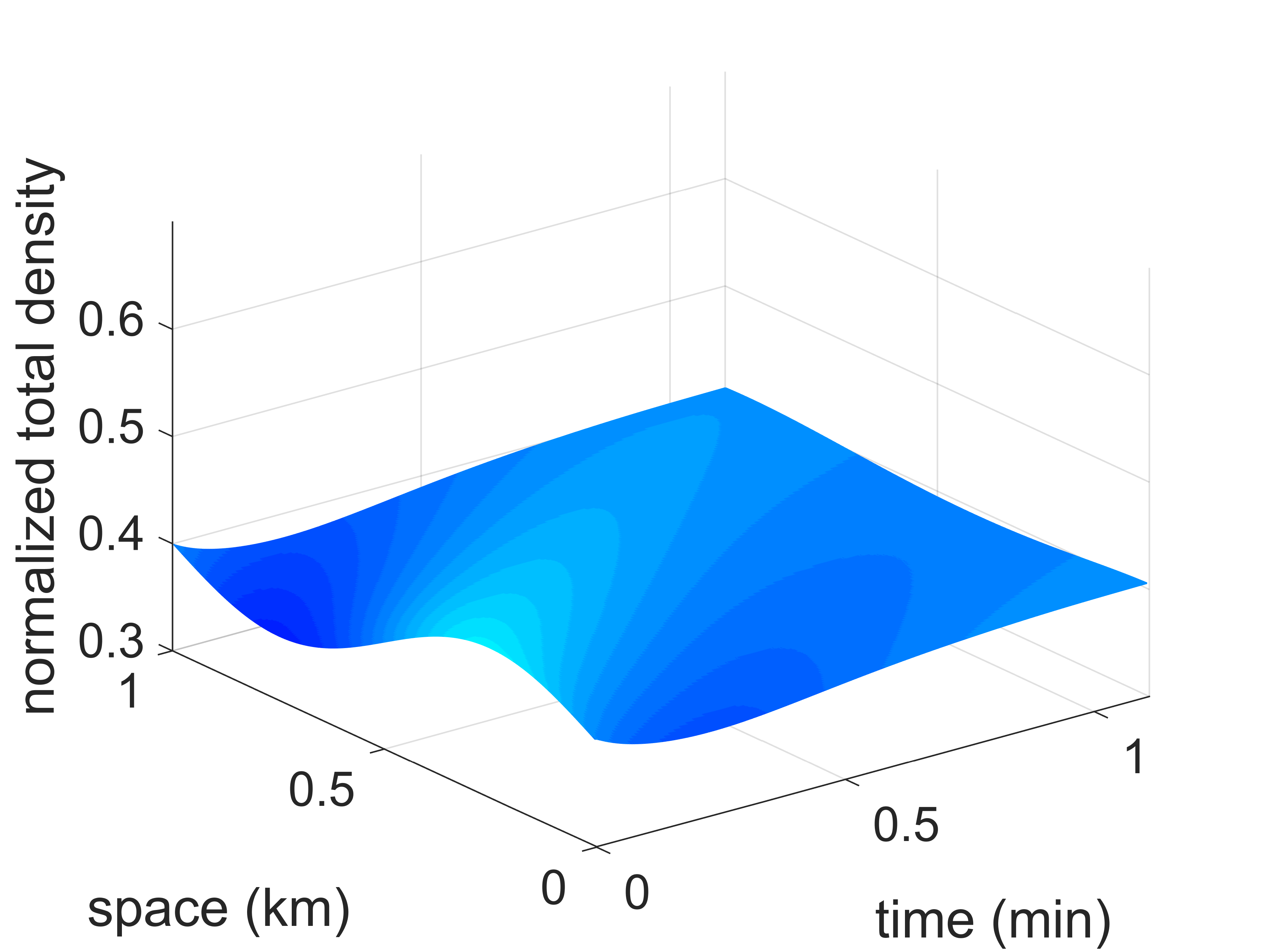}
	\end{minipage}
	\caption{Evolution of normalized total density when $\beta=0$, $\bar{\rho}^{\text{TOT}}=0.4\rho_{\text{jam}}$, 0\% AV (left) and 30\% AVs (right)}
	\label{fig:simu}
\end{figure}

Figure \ref{fig:simu} compares the total density evolution between a stable example and an unstable example. When the total density is $\bar{\rho}^{\text{TOT}}=0.4\rho_{\text{jam}}$, the pure HV traffic is unstable while 30\% AVs can stabilize the mixed traffic. In the former case, the initial perturbation on the total density grows up and develops a shock; In the latter case, the same initial perturbation decays and the total density converges to a uniform flow.

\begin{figure}[htbp]
	\centering
	\includegraphics[width=.18\textwidth]{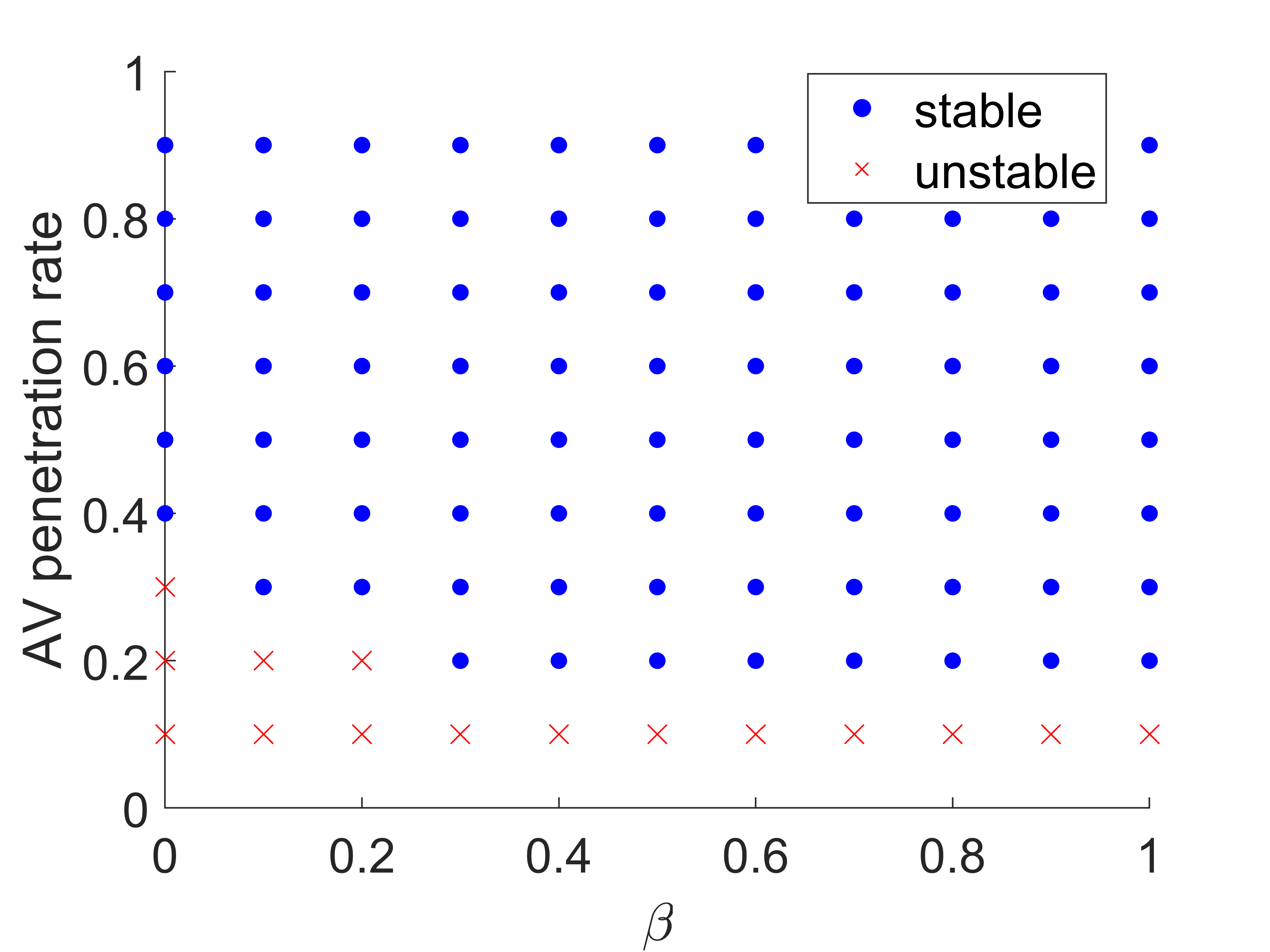}
	\caption{Stability region for the third group of experiments}
	\label{fig:beta}
\end{figure}

In the third group of experiments we fix the total density $\bar{\rho}^{\text{TOT}}=0.5\rho_{\text{jam}}$ and vary the AV's penetration rate and the parameter $\beta$. Then we plot the results in the phase diagram between $\beta$ and the AV's penetration rate, see Figure \ref{fig:beta}. We observe that for any fixed $\beta$, increasing AV's penetration rate makes the traffic more stable. When the AV's penetration rate is fixed but higher than 20\%, increasing $\beta$ makes the traffic more stable. This means that when AV is more sensitive to HV, the traffic becomes more stable.

\section{Conclusion}

This paper presents continuum traffic flow models for both pure AV traffic and mixed AV-HV traffic. The pure AV traffic is modeled by a mean field game and the linear stability analysis shows the traffic is always stable. The mixed AV-HV traffic is modeled by a coupled MFG-ARZ system. To demonstrate the mixed traffic stability analysis, three groups of numerical experiments are performed. In particular, we characterize the stability regions over AV density and HV density as well as over total density and AV's penetration rate in the mixed traffic. We also quantify the impact of the AV controller parameter on traffic stability. In future work, we plan to develop analytical stability analysis for mixed traffic and discuss the relation between more general AV controller designs and stability under different types of AV-HV interactions.


\section*{APPENDIX}
Apply Fourier analysis to \eqref{mfg_lin}, denote $\hat{\rho}(\xi,t)$ and $\hat{u}(\xi,t)$ the Fourier modes of $\tilde{\rho}(x,t)$ and $\tilde{u}(x,t)$, $\xi=\frac{2k\pi x}{L}$ ($k\in\mathbb{Z}$). For any frequency $\xi$:
\begin{align}
\begin{cases}
	\hat{\rho}_t+i\xi(1-\bar{\rho})\hat{\rho}+i\xi\bar{\rho}\hat{u}=0,\\
	\hat{u}_t+i\xi(\bar{\rho}-1)\hat{\rho}+i\xi(1-2\bar{\rho})\hat{u}=0.
\end{cases}\label{eq:ode}
\end{align}

It is an ODE system with the initial condition $\hat{\rho}(\xi,0)=\hat{\rho}_0(\xi)$ where $\hat{\rho}_0(\xi)$ is the Fourier transform of $\tilde{\rho}_0(x)$ and the terminal condition $\hat{\rho}(\xi,T)+\hat{u}(\xi,T)=0$. The linear PDE system \eqref{mfg_lin} is stable in $L^2$ norm if and only if there exists a universal constant $C>0$ such that for any $T>0$ and $\xi$, the solution of the ODE system \eqref{eq:ode} on $[0,T]$ satisfies:
\begin{align}
 	|\hat{\rho}(\xi,t)|^2+|\hat{u}(\xi,t)|^2\leq C|\hat{\rho}_0(\xi)|^2,\ \forall t\in[0,T]\label{eq:fin}.
\end{align} 

The ODE system \eqref{eq:ode} is homogeneous. We can assume without loss of generality that $\hat{\rho}_0(\xi)=1$.
To check the condition \eqref{eq:fin} we directly solve this boundary value problem of the ODE system \eqref{eq:ode}. Denote $r=\sqrt{\bar{\rho}  (5 \bar{\rho} -4)}$, $\eta=\xi t$, $\lambda =\xi T$ and $S=\exp \left(-\frac{1}{2} i \eta  \left(r-3 \bar{\rho} +2\right) \right)$,
the solution is:
\begin{align}
	\hat{\rho}(\xi,t)&=S\frac{(r+\bar{\rho})e^{i r  \eta}+(r-\bar{\rho})e^{i r  \lambda}}{r+\bar{\rho}+(r-\bar{\rho})e^{i r  \lambda}},\\
	\hat{u}(\xi,t)&=-S\frac{(r+3\bar{\rho}-2)e^{i r  \eta}+(r-3\bar{\rho}+2)e^{i r  \lambda}}{r+\bar{\rho}+(r-\bar{\rho})e^{i r  \lambda}},
\end{align}
when $\bar{\rho}\neq\frac45$ or $\hat{\rho}(\xi,t)=e^{\frac15i\eta}\frac{5i-2\eta+2\lambda}{5i+2\lambda}$ and $\hat{u}(\xi,t)=-e^{\frac15i\eta}\frac{5i-\eta+\lambda}{5i+2\lambda}$
when $\bar{\rho}=\frac45$.

Define:
\begin{align}
	E_{\bar{\rho}}(\lambda)=\max_{0\leq\eta\leq\lambda\text{ or }\lambda\leq \eta\leq 0}\left[|\hat{\rho}(\xi,t)|^2+|\hat{u}(\xi,t)|^2\right].
\end{align}
Then to check \eqref{eq:fin} it suffices to check the boundedness of the function $E_{\bar{\rho}}(\lambda)$ for all $0<\bar{\rho}<1$. We do this by computing the values of $E_{\bar{\rho}}(\lambda)$ from discrete values of $\bar{\rho}$ and $\lambda$. 
The computation shows that for any $\bar{\rho}$, $E_{\bar{\rho}}(\lambda)$ is bounded when $|\lambda|\to\infty$. 



\section*{ACKNOWLEDGMENT}
The authors would like to thank Data Science Institute from Columbia University for providing a seed grant for this research.

\bibliographystyle{IEEEtran}
\bibliography{IEEEabrv,survey}

\end{document}